# A Visual Sensitivity Analysis for Parameter-Augmented Ensembles of Curves

**Ribés, Alejandro**
EDF Lab Paris-Saclay, 7 Bd Gaspard Monge, 91120 Palaiseau, France
alejandro.ribes@edf.fr

**Pouderoux, Joachim**
Kitware, 6 Cours André Philip, 69100 Villeurbanne, France
joachim.pouderoux@kitware.com

**Iooss, Bertrand**
EDF Lab Chatou, 6 Quai Watier, 78401 Chatou, France
bertrand.iooss@edf.fr

**ABSTRACT**

*Engineers and computational scientists often study the behavior of their simulations by repeated solutions with variations in their parameters, which can be for instance boundary values or initial conditions. Through such simulation ensembles, uncertainty in a solution is studied as a function of the various input parameters. Solutions of numerical simulations are often temporal functions, spatial maps or spatio-temporal outputs. The usual way to deal with such complex outputs is to limit the analysis to several probes in the temporal/spatial domain. This leads to smaller and more tractable ensembles of functional outputs (curves) with their associated input parameters: augmented ensembles of curves. This article describes a system for the interactive exploration and analysis of such augmented ensembles. Descriptive statistics on the functional outputs are performed by Principal Component Analysis projection, kernel density estimation and the computation of High Density Regions. This makes possible the calculation of functional quantiles and outliers. Brushing and linking the elements of the system allows in-depth analysis of the ensemble. The system allows for functional descriptive statistics, cluster detection and finally for the realization of a visual sensitivity analysis via cobweb plots. We present two synthetic examples and then validate our approach in an industrial use-case concerning a marine current study using a hydraulic solver.*

# INTRODUCTION

In this article, simulation refers to the application of computational models to the study and prediction of physical events or the behavior of engineered systems. In this context, the modern usage of simulation tools has improved and grown to a point that has far exceeded many expectations. That remarkable change has come about mainly because of developments in the computational sciences and the rapid advances in computing equipment. Computer models help engineers to forecast the behavior of the system under investigation in conditions that cannot be reproduced in physical experiments (e.g. accidental scenarios), or when physical experiments are theoretically possible but at a very high cost. To improve and have a better hold on these tools, it is crucial to be able to analyze them under the scopes of sensitivity and uncertainty analysis [1-3]. In particular, sensitivity analysis aims at identifying the most influential parameters for a given output of the computer model and at evaluating the effect of uncertainty in each uncertain input variable on model output [4, 5].

A probabilistic uncertainty study consists of evaluating the computer model on a large size statistical sample of model inputs (which follow a joint probability distribution), then analyzing all the results (the model outputs) with specific statistical tools. The result of such a family of runs is called *ensemble*, and each individual run is called a *member*. Ensembles are multivariate, which means that a simulation is run several times with varying parameters. Their members are multidimensional (both in space and time) and multivalued (several quantities such as temperature, pressure or velocity are considered). The usual way to deal with these kinds of outputs (as a temporal function, a spatial map or a spatio-temporal output) is to limit the analysis to several probes in the temporal/spatial domain [1-4]. To deal with this problem, taking ideas from the visualization community seems particularly interesting. Indeed, one of its current challenges is how to deal with the multivariate nature of the ensembles [6-8]. Furthermore uncertainty visualization has been long advocated as one of the top challenges in visualization [9-11].

The goal of this work is to propose methodologies and tools for researchers and engineers performing uncertainty studies by analyzing ensembles. An example of such a strategy can be a hydraulics engineer studying results generated by a multi-run finite-elements simulation. In this case, the ensemble could be a fixed 3D mesh for all members and a varying field (temperature, water height, pressure, etc.) that depends on the experimental design used to sample the parameters controlling the simulations. Thus, when the engineer applies a probe on a node of the mesh she/he obtains not the evolution of a quantity (temperature, water height, pressure, etc.) over time but another smaller ensemble of functional outputs or curves. We should then not only deal with an ensemble of functional outputs but also with their associated simulation parameters. We call this kind of data an *augmented ensemble of curves*.

A non-augmented ensemble of curves presents already a first problem of visual clutter, which is well known in the visualization community. When a large number of curves are superposed to one another, the overall perception of the graphs is lost, and the user cannot analyze the ensemble. As an example, Figure 1 depicts 1,500 curves coming from different runs of the same numerical simulation (from a hydraulics application). When looking at the overall behavior of an ensemble of curves, such as Figure 1, the first set of basic questions that arise are the following:

- What is the median curve?
- Can we define some confidence interval curves containing most of the curves (as done usually for scalar random variables with the boxplot tool)?
- Can we detect some abnormal curves, in the sense of a strong difference from the majority of the curves (as outliers for scalar variables)?
- Are there some clusters, which correspond to different behaviors of the physical model that generated these outputs?

These questions can be answered by methods found in the recent technical literature by the way of Principal Component Analysis (PCA) methods, with a statistical viewpoint [12-14] or with a visualization viewpoint [15-17]. However, for augmented ensembles new challenges arise because a member of such ensemble consists of a set of input parameters (which drove a numerical simulation) and its associated functional output. First, the interactive exploration needs a methodology able to visually provide, to the analyst, the statistical structure of the curves and the identification of clusters. Second, if the clusters of functional outputs correspond to groups of coherent behaviors of the simulations, is it possible to visually study the relationship between these behaviors and the input parameters? This question implies the realization of a *visual sensitivity analysis* that we realize by linking the classical tool of the cobweb plot in sensitivity analysis [4] and the ensemble of curves visualization described before.

The following section lists the important and main previous works on the subjects covered by this paper. The third section explains the method used for estimating functional quantiles while the fourth section describes how to perform the visual sensitivity analysis. In the fifth section, applications of the methodology are given on toy examples and an industrial example. The two last sections provide a discussion on software implementation and a conclusion.

**BACKGROUND AND RELATED WORK**

First of all, our work relates to uncertainty and **sensitivity analysis**. In particular, global sensitivity analysis is an ensemble of techniques which aim to identify the influential and non-influential inputs on some computer model outputs [4]. In particular, quantitative global sensitivity analysis deals with a probabilistic representation of the input parameters to consider their overall variation range. Variance-based sensitivity measures, also called Sobol' indices [18], are currently the most popular method for global sensitivity analysis [5]. The principle of Sobol' indices is to decompose the variance of the output, $Y$, of the simulation into fractions, which can be attributed to each of the random model input $X_i$ (with $i=1,…,p$ where $p$ is the number of inputs). When $Y$ is a scalar output, these percentages are directly interpreted as measures of sensitivity. However, sensitivity analysis for large scale numerical systems that simulate complex spatial and temporal evolutions remains very challenging because of the treatment of uncertainty [2], the treatment of the functional nature of the output [19,20] and the large volumes of data that could be produced [21]. Our main contribution is the realization of *a visual sensitivity analysis*, linking the cobweb plot (a classical graphical tool in sensitivity analysis [4]) and the ensemble of curves visualization.

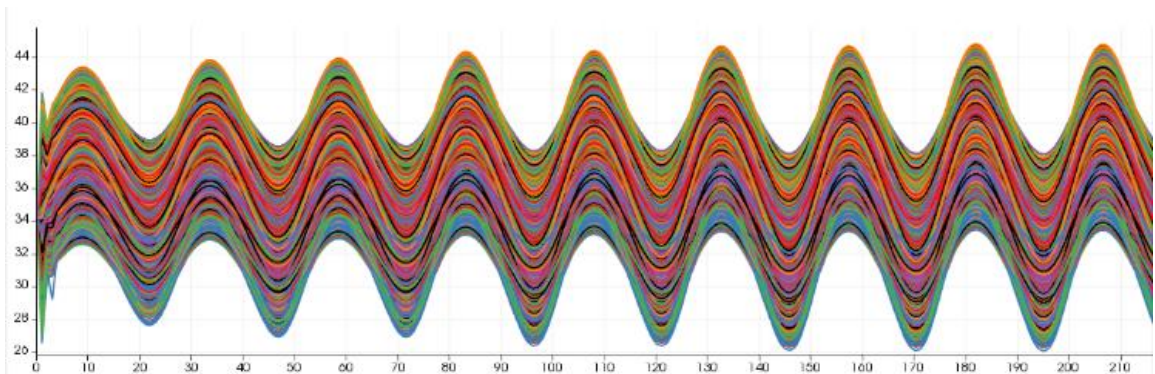

Figure 1. Raw visualization of curves coming from a multi-run hydraulics simulation: 1500 curves of water height evolving over time.

One of the difficulty when visualizing several one-dimensional curves is to avoid **visual clutter**. An interesting solution is given by [22] which proposes a "curve density estimation" directly in the curves' space. Other visualizations may represent overlaid function graphs as envelopes [23], semi-transparent graphs [24], and offer brushing techniques to highlight selected subsets of the functions [23-25]. Another approach offers a re-orderable matrix of time series charts [26]. The way these methods deal with visual clutter is very different from our approach, which offers quantified statistical information by calculating quantile bands and outliers.

In our work, we use the extension of the classical boxplot to functions: the **functional boxplot** proposed in [12,27]. A boxplot for scalar variables allows summarizing the main information of a data sample: median, first and third quartiles, and an interquantile-based interval which define the limit of non-outliers data. First step to build such boxplot is to rank data thanks to a statistical order or *data depth*; such order has first to be defined for functional data [28], which has led to numerous research works in the literature. The concept of functional data depth has been generalized to contours by [15], which displays boxplots for two-dimensional simulation data in weather forecasting and computational fluid dynamics. The so-called band depth, defined by [28], is particularly relevant for the goals of [15]. Band depth is defined on an ensemble of functions, the band depth of each function is the probability that the function lies within the band defined by a random selection of other functions from the distribution. The band depth is computed for each member of the ensemble and can be used, as described in [27], to visualize summary statistics for an ensemble of functions. Our methodology differs from [15] because we do not use data depth, the functional summary statistics are calculated through and alternative method based on a PCA projection. Furthermore, our focus is in augmented ensembles and not in ensembles of contours.

In our method, the functional curves are handled by reducing their dimension via projection, using PCA as in [12] (which is limited to the two first components). [13] has introduced this PCA-based approach for visualizing (but non-interactively) functional outputs of computer experiments. Later, [14] has extended the technique to selecting and modeling more than two PCA components by advanced statistical techniques. Our choice of using PCA is firstly motivated by our idea to jointly interact with the PCA-plane (defined in section "**Projecting on the PCA bivariate plane**") in which each function is represented by a 2d-point. Thus, in this aspect the method relates to dimensionality reduction techniques but using a human-in-the-loop approach; see [29] for a structured literature review and references on this field. Furthermore, the PCA-plane allows, at the same time, to calculate functional quantiles and to study the multimodal nature of *ensembles*. We remark that data depth based techniques of [27] do not deal with multimodality.

The PCA technique reduces the data dimension via a linear transformation. In some cases, such a transformation does not work due to the underlying structure of the data (see an engineering example in [30]). Non-linear dimension reduction techniques (non-linear PCA, kernel PCA, Riemannian manifold learning, locally linear embedding, etc., see [31]) can then be used with a certain increase in complexity and computational cost. The pragmatic approach consists in first applying a PCA, and then turning to non-linear methods if the data variability is not well captures by a small number of PCA components. Using non-linear dimension reduction techniques is beyond the scope of this paper and will be studied in future works.

Reference [17] presents a method for computing streamline variability plots. It consists on the transformation, via PCA, of the set of streamlines, into a lower dimension space in which clustering can be performed. Clustering by fitting geometric medians and confidence ellipses is performed in PCA- space. Finally, medians and ellipses are transformed back to the domain space and yield the variability plot of the streamline ensemble. Reference [16], of the same authors, applies similar techniques to ensembles of iso-surfaces. Our methodology also uses PCA in order to work into a lower dimension space and then re-project to the original dual space. However, the operations performed in the PCA-space are very different, [16, 17] perform clustering while we calculate HDR (High Density Regions). This operation necessitates an estimate of the empirical density function in the PCA space, which uses kernels thus avoiding fitting a parametric model (such as ellipses). HDR presents a unique and strict mathematical definition, and are at the core of the method to calculate quantitative and non-parametric variability of the data. HDR can also be used for clustering, in our current implementation they assist the users in this task.

**Brushing and linking** is extensively used in our system but we have no new contribution in this area, the interactions we use were already described in classical works such as [32,33]. Other works have applied these classical techniques to ensembles of functions, for instance [34] for the investigation of families of data surfaces, [35] to analyze 2D function ensembles in the development process of powertrain systems and [36] for the interactive visual exploration of large 3D scalar ensembles. These references demonstrate the necessity for a flexible visual analysis system that integrates many different linked views for making sense of this complex data. In this context, using brushing and linking and statistical aggregations, as [34,36] is seducing; we differentiate from these works because we do not perform any statistical aggregation, i.e. the computation of statistical moments, and prefer quantile analysis that does not introduce any hypothesis about the underlying data distribution.

We finally remark that our work is inscribed in the field of **visual parameter space analysis**. This approach was used by [37] in an interactive system called HyperMoVal that was designed to support model validation. These models related to the development of car engines for tasks which require a prediction of results in real-time. Other examples of this visual analysis include [38], which combined a sensitivity analysis with a linked multi-dimensional visualization providing a way to analyze the behavior of an artificial neural network. Reference [39] addresses the problem of parameter-finding in image segmentation by visually guiding the user towards areas that need refinement, in a sparse sampled parameter space, by placing additional sample points. In a second stage the user navigates through the parameter space in order to determine areas where the response value (goodness of segmentation) is high.

Reference [40] presents a conceptual framework in which six typical analysis task can be performed: optimization, partitioning, fitting, outliers, uncertainty and sensitivity. Numerous examples exist of space analysis for optimization, such as [37] or [39]. Our work differs from most references [34-38] in the analysis tasks we focus on (detecting outliers, partitioning and visualizing sensitivities). In fact, our main contribution is the realization of a visual sensitivity study in the context of time-evolving numerical simulations.

## ESTIMATING FUNCTIONAL QUANTILES AND OUTLIERS

The estimation of the quantiles and outliers of an ensemble of functions is performed on a plane defined by the first two vectors of its Principal Component Analysis (PCA). The method is divided into 3 main steps:
1. Project the functions into the PCA bivariate plane;
2. Perform the estimation of the Probability Density Function (PDF) on this plane, which allows for the estimation of Highest Density Regions (HDR) which boundaries are isoprobability contours;
3. Project the HDR boundaries back into the space of curves. The functional quantiles and outliers are then computed.

### Projecting on the PCA bivariate plane

PCA is a technique of dimensionality reduction, whose purpose is to represent the source data into a new space of lower dimensions. It is mathematically defined as an orthogonal linear transformation, which maps the data to a new coordinate system such that the greatest variance comes to lie on the first coordinate (called the first principal component), the second greatest variance on the second coordinate (called the second principal component), and so on. Projecting the curves into the PCA bivariate plane means that only the first and second principal components are kept, thus each curve is represented by a point in a two-dimensional space (the bivariate plane).

There is an underlying condition for PCA reduction: the transformation should keep enough information about the source data, while allowing to simplify the analysis. In this article, we limit our scope to a bivariate plane for simplicity but an extension to larger dimensions is possible. Moreover, our examples have a high explained variance using the first two PCA components (which means that the two-dimensional reduced PCA basis correctly reproduces the overall variability of the ensemble of curves). This explained variance $v_e$ is calculated from the $n$ singular values of the covariance matrix $(\sigma_1, \sigma_2..., \sigma_n)$ in the following way:

$$v_e = \frac{\sigma_1 + \sigma_2}{\sigma_1 + \sigma_2 + \ldots + \sigma_n}$$

where $\sigma_1 \geq \sigma_2 \geq \ldots \geq \sigma_n$ are the ordered singular values representing the importance of the variance of each principal component. In our implementations, this quantity is systematically visualized on the diagrams containing the bivariate plane.

As said before, working with three or more PCA components is also possible. We have prepared a prototype that uses an interactive Scatter Plot Matrix view, based on [41], for the interaction with the *n*-dimensional space of PCA dimensions. However, the real challenge in this case is the computational complexity associated to the estimation of densities in higher dimensions. The work of [14] allows to solve this problem.

**Highest Density Regions method**

Once the functional variable has been transformed to a fewer component space (in our case the bivariate plane), the next step is to estimate the quantiles, which also allows outliers to be detected. Conceptually, two basic operations should be performed to build the quantiles: i) create a density map on the plane and ii) calculate iso-probability curves of this map. So as to implement these operations we follow the method of Highest Density Regions of [12]. Full mathematical details are given in [12] but the principle of this method is to assimilate observations in the bivariate plane of principal components to the realizations of a random vector with density *f*. By calculating an estimate of the density map *f*, the quantiles can then be computed.

We start from the sample $(X_i)_{i=1\ldots n}$ which stands for *n* observations of the vector *X* of dimension *p*=2. The following smoothing process is used:

$$\hat{f}(X) = \frac{1}{n}\sum_{i=1}^{n} K_H(X - X_i)$$

where $K_H$ is the Gaussian smoothing kernel which writes

$$K_H(X) = |H|^{-1/2} K(H^{-1/2} X)$$

with $K(X) = \frac{1}{2\pi} exp\left(-\frac{1}{2}\langle X, X \rangle\right)$ the "standard" Gaussian kernel and H the matrix containing the smoothing parameters (extension in *p* dimensions of a smoothing parameter *h* in dimension 1). Depending on this matrix (diagonal or not), some preferential smoothing directions can be chosen. In our implementation, we first generate a grid covering the bivariate plane, which is initialized to 100x100 but can be customized by means of the user interface. We subsequently apply an isotropic kernel, whose width is automatically initialized by use of the rule of Silverman [42].

Once the estimate of the density map $\hat{f}$ is obtained, the Highest Density Regions method (HDR) gives a description of important statistical information. It is defined as

$$R_\alpha = \{Z : \hat{f}(Z) \geq f_\alpha\}$$

where $\alpha$ is the order of the quantile we choose and the quantile value $f_\alpha$ is such that $\int_{R_\alpha} \hat{f}(Z)dZ = 1 - \alpha$, that defines the region with probability coverage $1 - \alpha$. It means that $f_\alpha$ is such that all points within the region $R_\alpha$ have a higher density estimate than any of the points outside the region, hence the name *highest density region*. For a density map, the HDRs can be considered as regions bounded by contours, with an expanding coverage as α decreases. In our implementation we compute two HDRs:

- An inner HDR with probability coverage α = 50% that corresponds to the central interquartile zone;
- An outer HDR where α can interactively be modified via a slider (default value α = 95%).

We consider that all points excluded from the outer HDR can be some outliers. This allows to interactively changing the threshold that is used to compute the outliers.

**Back into the curves space**

Once both HDRs are calculated, we would like to see them in the original space of the curves. We recall that a point in the PCA plane corresponds to a curve. However, HDRs represent areas of the PCA plane, which boundaries are contours. It is then necessary to run an algorithm that converts these contours into their associated functional quantiles. We propose a first exact algorithm by traversing all points corresponding to the discretized boundary and choosing the maximum and minimum values on the curves space. This process generates functional quantiles that are not necessarily existing curves on the ensemble.

Figure 2 shows the analysis of the dataset shown in Figure 1: the 50% inter-quantile area is represented in light red, dark red is used for the 95% inter-quantile zone. The median curve (in black) is calculated by finding the point in the bivariate plane that presents the highest value of the density map.

**Information contained in the bivariate plane**

The bivariate PCA plane presents some important characteristics worth discussing:
- It provides a data reduction and visually understandable representation, where each curve is represented by a point;
- The probability density map associated to the bivariate plane allows for the calculation of the median curve, functional quantiles and outliers;
- The density map conveys important information about the modality of the curves dataset. Indeed, statistical multimodality is normally associated with a mixture of unimodal distributions. Each of the underlying modes defines different behaviors of the curves and thus the original data can be divided in clusters. The HDR exposes the mono or multi-modal nature of the dataset. If an HDR is formed by disjoint areas, the distribution is multimodal.
- 

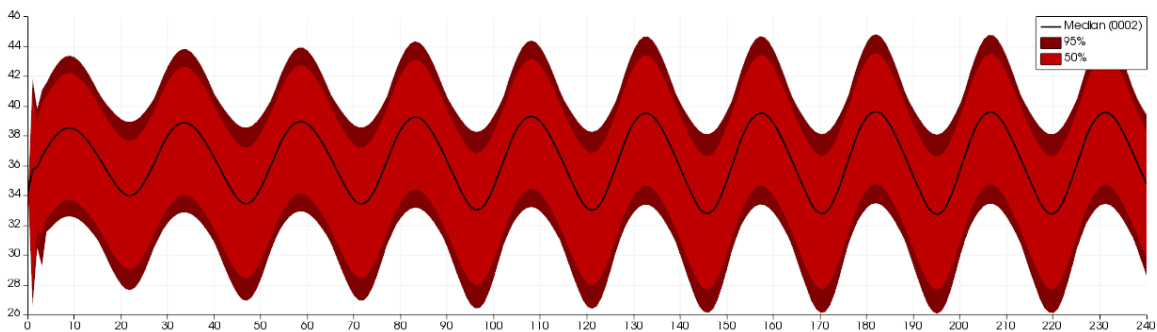

Figure 2. The functional quantiles and the median (black line) of the multi-run hydraulics simulation curves shown in Figure 1.

In the presented method, both the diagram containing the curves or *functional boxplot* and the bivariate PCA plane are jointly visualized by use of a brushing and linking strategy. In our system brushing corresponds to a selection operation. Thus, we offer tools to select individual or subsets of points in the bivariate plane, which highlights its corresponding curves. The opposite schema "select curve, highlight point" is also available. We also set meaningful limits to the exploration by drawing two iso-probability curves. First contour is fixed at probability 50% while the outer one is controlled by a slider in the user's interface. We finally define a blue vs red colormap to help the visual interaction with the plane (blue meaning low probability and red high probability).

**LINKING WITH THE INPUT PARAMETERS**

The proposed methodology also allows for the study of the *augmented ensemble*, each member of such *ensemble* is a couple consisting of the input parameters and a functional output. Thus, a member is represented as a couple ($p_N$, $f$) where $p_N$ is a list of $N$ parameters and $f$ is a function. We will constrain our current study to one-dimensional $f$ functions (curves). We also remark that the number $N$ of input parameters is not necessarily small. Current numerical simulations can easily present $N$=50. The whole ensemble is represented as ($p_N$, $f$)$_M$ where $M$ is the number of couples or, equivalently, the number of members in the ensemble. We remark that input parameters and functional outputs both possess a common index $M$ on the ensemble. Thus, it is technically possible to share the same selection strategy between them. Indeed, our linking strategy allows the use of fully coupled diagrams in order to interact, simultaneously, with the input parameters and functional outputs of the multi-run simulations.

One important consequence of this joint visualization is that it allows for what we define as a *visual sensitivity analysis*. As a matter of fact, multi-run simulations are often used to determine the impact of input parameters on the results of the simulations, which is called a sensitivity analysis [4]. Our system propagates the selections performed on the bivariate plane or on the functional boxplot to the diagrams associated to the input parameters. This allows the exploration of complex input-output relationships.

**RESULTS**

In this section, we discuss experimental results to demonstrate the utility of the proposed PCA-based functional boxplot, which allows:
1. To study the variance of the curves generated from a multi-run numerical simulation;
2. To detect functional outliers;
3. To identify clusters of curves that correspond to different behaviors of members of the ensemble;
4. To perform a **visual sensitivity study**.

The discussion is started with two synthetic examples and then an industrial use-case concerning a marine current study using a hydraulic solver is presented.

**Oscillating tangents**

Our first synthetic example consists of an ensemble of time-oscillating analytical functions coming from the following equation:

$$y(t) = atan(X_1)cos(t) + atan(X_2)sin(t)$$

where $X_1$ and $X_2$ are the input parameters and t represents the time, which is regularly sampled in the interval [0, 2π]. Thus we generate an ensemble of 400 curves by Monte Carlo sampling of both $X_1$ and $X_2$ based on a uniform distribution in the interval [-7, 7].

In Figure 3, we show some results concerning this ensemble of temporal oscillating functions. On the top panel (a) all 400 generated curves are shown. Figure 4 (b) shows the result of a user interaction with the PCA bivariate plane of the 400 curves, where a blue to red colormap is applied. The explained variance is equal to 100%. This surprising result is explained by the fact that curves of the oscillating tangents function are regular sinusoids only tuned by their amplitude and frequency.

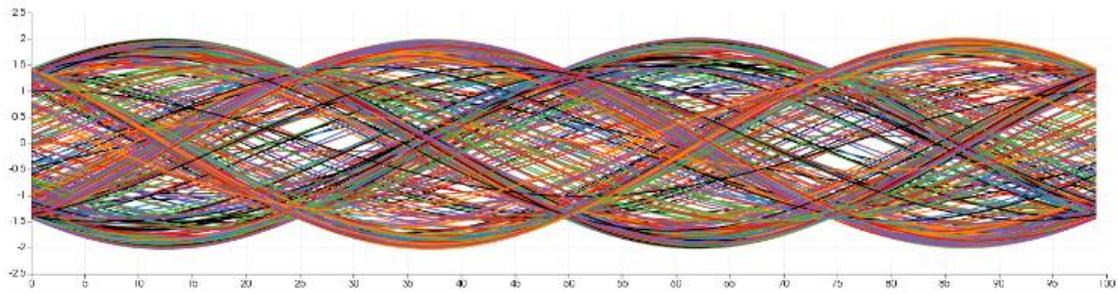

(a)

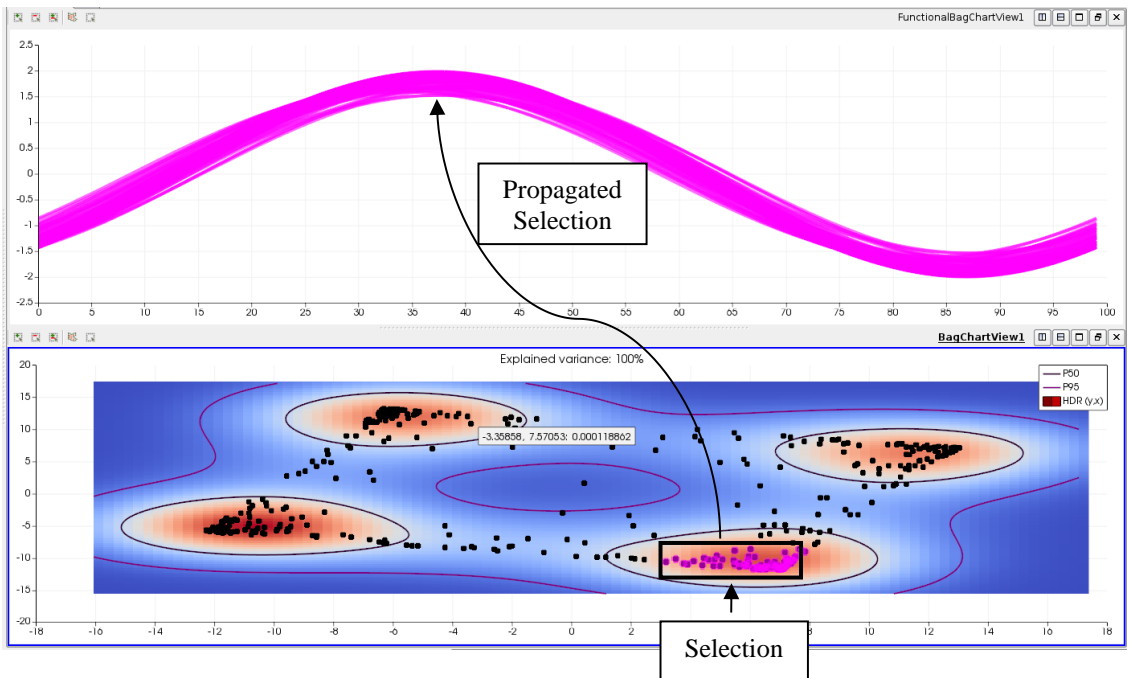

(b)

Figure 3. Top panel **(a)** shows 400 curves generated by the temporal oscillating tangents experiment. Bottom panel **(b)** shows the results of a user interaction where the analyst has selected (pink points) one of the clusters of the PCA plane.

We can see that four clusters appear, indicating a multi-modal structure of the oscillating curves. The analyst has selected one of these cluster, then the propagated selection on the curves is highlighted, this selection corresponds to variations of the same oscillating mode. This example demonstrates the interest of visualizing and interacting with the PCA bivariate plane in the context of a partitioning task [40]. Understanding a multimodal ensemble of curves is indeed a complex analysis task.

### Campbell 1D functions

Our second synthetic example is inspired by [43, 19]. It consists of an ensemble of analytical functions that evolve in time. This dataset is generated by use of the following equation:

$$y(\tau) = 10 + X_1 exp\left(\frac{-(\tau - 10X_2)^2}{k_1 X_1^2 + X_3^2}\right) + X_2 X_4 exp(k_2 X_1 \tau)$$

where $X_1$, $X_2$, $X_3$ and $X_4$ the input parameters and $\tau$ is a one by one regularly sampled variable in the interval [-90, 90]. The quantities $k_1$ and $k_2$ are constant, fixed to 60 and 0.002 respectively. [43] has introduced a slightly different version of this function in order to test simple sensitivity analysis tools when model outputs are 1D curves (understanding the role of each of the four inputs on the translation from left to right of the curve, on the shape of the curve peak and on the curve tail behavior). From this, [19] has calibrated a function (called Campbell2D) in order illustrate tools of sensitivity analysis when model outputs are 2D spatial functions (with strong spatial heterogeneities, sharp boundaries, and very different spatial distributions of the output values according to the **X** values).

We generate an ensemble of 100 curves by Monte Carlo sampling based on a uniform distribution in the interval [-1, 5], the same sampling is used for all $X_i$. In the upper panel (a) of Figure 4, we show all 400 curves generated by the Monte Carlo sampling. At time 80 an event occurs and part of the curves diverge from its original tendency while the others keep with its original behavior. This can be easily understood by looking at the functional interquantile areas and to the median curve, presented in the bottom panel (b) of Figure 4. Indeed, the median curve and the 50% interquantile area are not modified by the event while the upper limit of the 95% interquantile area rises up. By looking at this representation, an analyst avoids visual clutter and easily understands that the event at time 80 affected only the evolution in time of the top 25% of the curves. The explained variance by the two PCA components is equal to 97%.

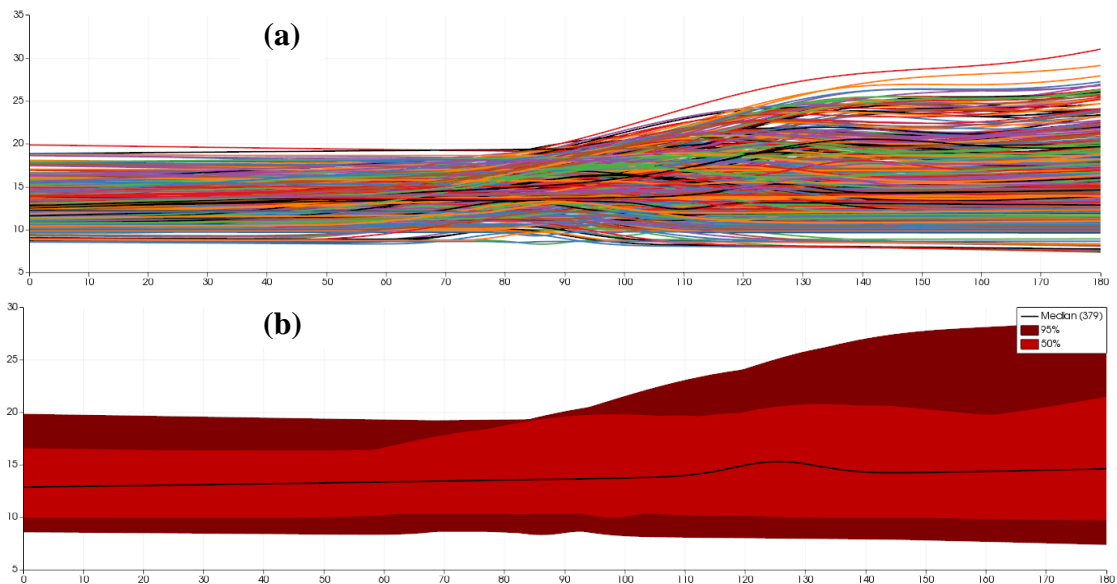

Figure 4. Top panel **(a)** shows 400 curves generated by the modified 1D Campbell function experiment. Bottom panel **(b)** shows the corresponding median curve and interquantile areas at 50% and 95% probability.

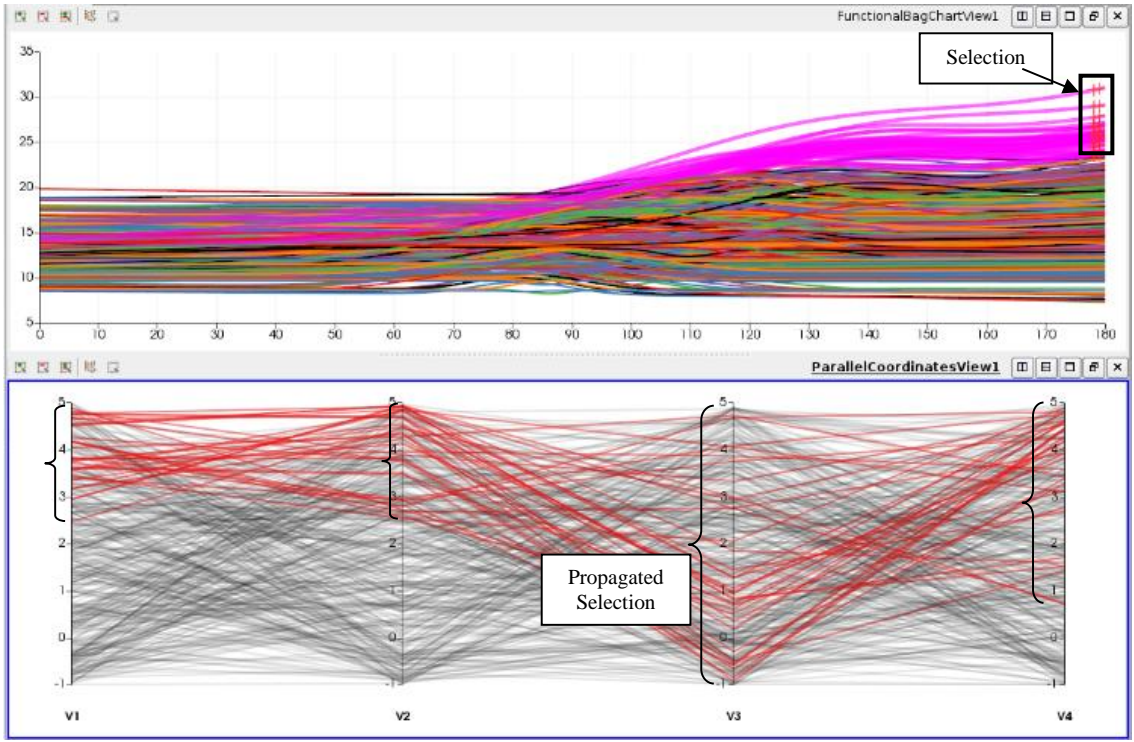

Figure 5. Realization of a visual sensitivity study over a synthetically generated ensemble using the modified 1D Campbell function. Two interactive linked diagrams are presented: the diagram on the top contains the analysis of the outputs (the curves), while the bottom parallel coordinates diagram represents the four input parameters of the Campbell function. The analyst has selected a group of curves on the top diagram thus this selection is propagated. We superpose "left bracket" symbols to the parallel coordinates' diagram in order to visually reinforce the dispersion of the propagated selection, which is a measure of sensitivity.

Once stated that there is a specific group of temporal evolving functions which behavior is modified by the event, the analyst is interested in knowing if some of the input parameters are responsible for this behavior. It is then possible to perform a *visual sensitivity study* by selecting, in the functional boxplot diagrams, all the curves ending in the upper part of the 95% interquantile area. Then the system propagates the selection to the diagrams dealing with the input parameters. The result of this operation is shown in Figure 5. In this figure, two interactive linked diagrams are presented, the diagram on the top contains the analysis of the outputs while the bottom parallel coordinates' diagram represents the inputs. The interpretation of the diagrams is straightforward. Indeed, it is possible to visually assess the importance of each parameter by looking at the axis of the parallel coordinates' diagram. In this case, $X_1$ and $X_2$ presents a high degree of concentration of the selection, thus they strongly influence the outputs. Using this simple criterion of "visual dispersion" the parameters can be ordered by importance ($X_1$, $X_2$, $X_4$, $X_3$), which is one of the main objectives of a sensitivity analysis.

In our system, the analyst "asks a question" by using the selection. In this case the question was: "which parameters generated the curves ending with the highest values?". Of course, numerous other questions are possible, based on the selection on the Functional Boxplot, the bivariate plane or other diagrams linked to the ensemble data. This example demonstrates that our system can perform a visual sensitivity analysis. This kind of fast and informative exploration could be performed before a formal sensitivity analysis, such as the computation of Sobol' indices (see [5] for this methodological point of view).

**A hydraulics study-case**

Our study-case concerns a maritime model of Alderney Ray (or Raz Blanchard in French), which is a strait that runs between Alderney (UK) and Cap de la Hague (France), a cape at the northwestern tip of the Cotentin peninsula in Normandy. This strait presents one of the fastest marine currents in Europe; the current is intermittent, varying with the tide, and can run up to about twelve knots during equinoctial tides.

A study was performed in order to calibrate a hydrodynamic model, which is typically an engaged and difficult process due to the complexity of the flows and their interaction with the shoreline, the seabed and the islands. Thus, it was essential to understand the relationship between the modelling calibration parameters and the simulated state variables which are compared to the observations. A sensitivity analysis using Sobol' indices was a necessary step prior to calibration. In this context, several multi-run studies were performed. In this section we focus on a particular 1,500 runs study where five parameters were varied:

- Two coefficients of friction (CF1 and CF2) modeling the interaction with the seabed.
- One "SeaLevel" representing the vertical distance from the surface to the seabed.
- Two parameters for tidal modeling: the tidal range (vertical variation range) and the tidal velocity.

The maritime model includes Alderney and the tip of the Cotentin peninsula and covers an area roughly 55 km x 35 km. The finite element mesh is composed of 17,983 nodes and 35,361 triangular elements. The mesh size varies from 100 m, at the shoreline and within the areas of interest, to 1.8 km offshore (western and northern sectors of the model). The computations were performed by the open source fluid dynamics solver TELEMAC [44] (http://www.opentelemac.org/) that generated fields such as velocity, pressure, and water height. We extracted 1,500 curves of this multi-run study by use of a probe in one of the nodes of the mesh; this leads to the curves shown in Figure 1.

Figure 6 shows the result of two interactions. A functional boxplot containing the analysis of the 1,500 curves is linked to the input parameters that are represented in a parallel coordinates' diagram. In this figure, the analyst explores the relationship between the functional outputs and the parameter "Sea Level". On the top panel (a) of Figure 6, the analyst selects the highest values of "Sea Level" while on the bottom panel (b) the lowest values are selected. By looking at the propagated selections on the functional boxplots (in orange), it is easily understood that "Sea Level" behaves like a vertical offset on the oscillating curves generated by the tide. The analyst thus understands that "Sea Level" strongly influences the simulations results. This is coherent with the formal sensitivity analysis that was also performed. Sobol' indices were computed and they show that the parameter "Sea Level" strongly influences (around 97%) the outputs while the others present little influence. In addition, the information shown in Figure 6 is richer than the scalar Sobol' indices. Sobol' indices reveal the strong influence of the parameter "Sea Level" while Figure 6 underlines the way that this influence is performed (by applying a vertical shift to the tide).

Hydraulics engineers were also interested in using our system to study or verify which parameters do not influence the functional outputs. This step is fundamental for model reduction where a parameter is taken out of a model when it is considered as non-influential. Figure 7 shows the result of selecting the highest values of CF1 (one of the coefficients of friction of the seabed). We observe that its propagated selection on the functional boxplot is visually disperse, which indicates that CF1 has no influence in the behavior of the outputs. This again is coherent in respect to the Sobol' indices-based sensitivity analysis. Moreover, physicist performing the study confirmed that CF1 and CF2 should be non-influential in this case because the seabed is too deep for its friction to have an effect on the sea surface. The explained variance by the two PCA components is equal to 99%.

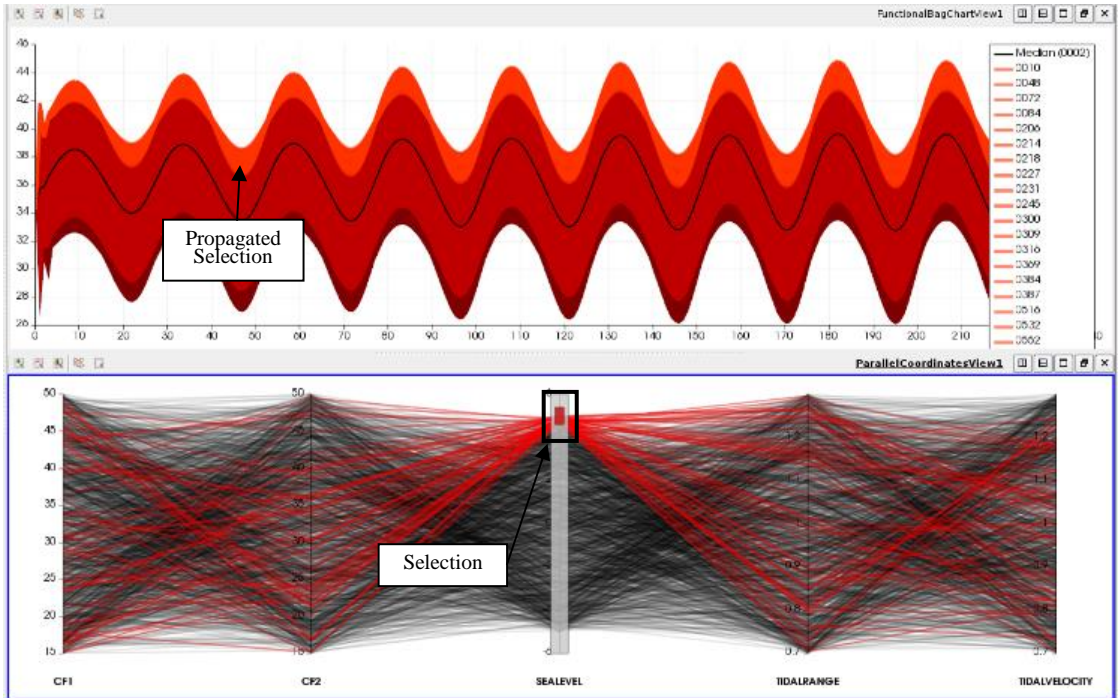

**(a)** High values of "Sea Level"

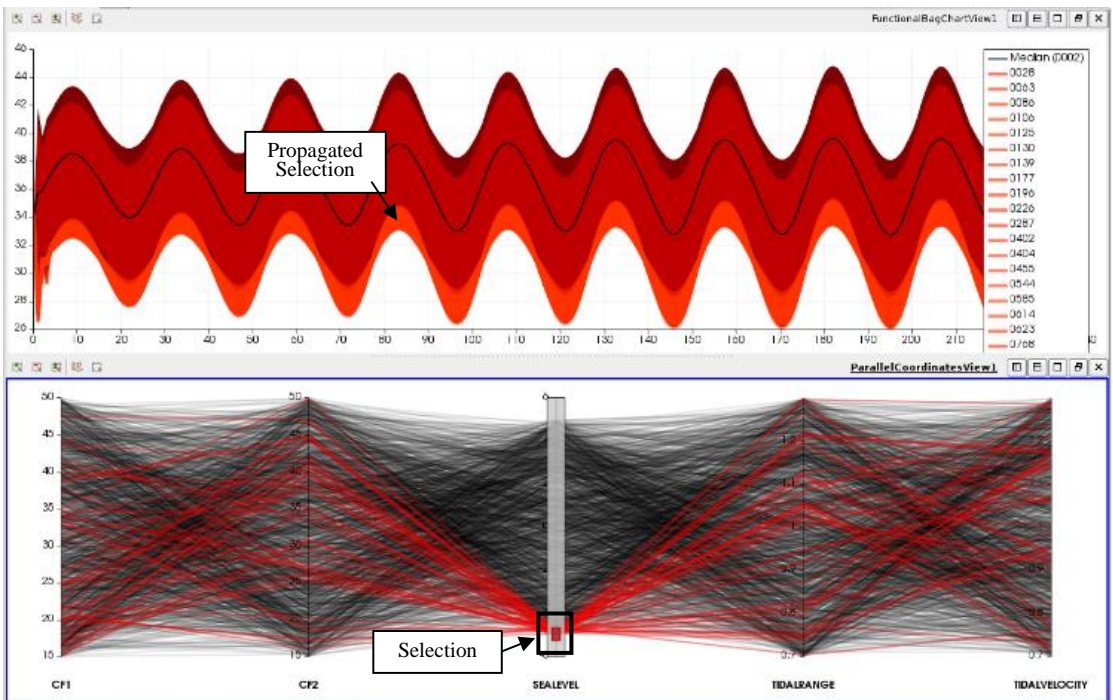

**(b)** Low values of "Sea Level"

Figure 6. Interactive exploration of the relationship between the functional outputs and the parameter "Sea Level" in a marine hydraulics multi-run study, which shows that this parameter applies a vertical shift to the tide.

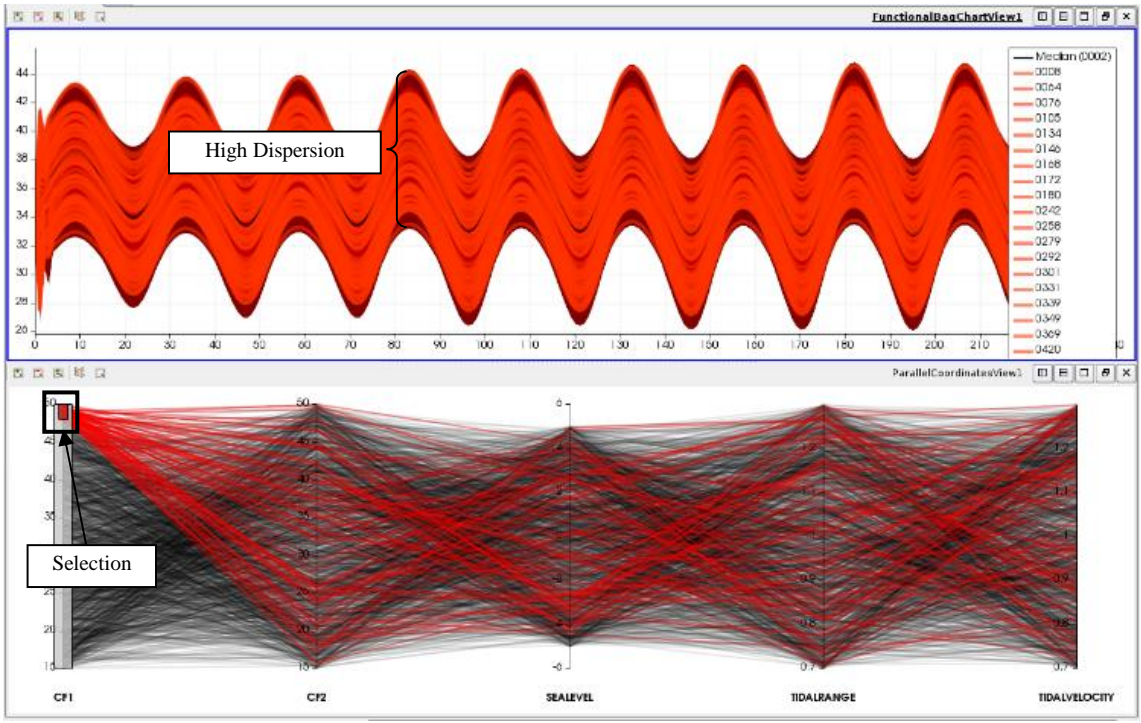

Figure 7. Interactive exploration of the relationship between the functional outputs and "CF1" (coefficient of friction 1) shows that this parameter is not relevant for this study, due to of the high dispersion of the propagated selection.

Finally, Figure 8 shows a more subtle result. The analyst interacted with the propagated selected curves of the bottom panel (b) of Figure 6. Our system allows refining selections, then sub-ensembles of the curves with low "Sea Level" values were selected and a second order or indirect effect was observed. Figure 8 illustrates this second order effect by selecting:
(a) low "Sea Level" and high "Tidal Range" values,
(b) low "Sea Level" and low "Tidal Range" values.
Comparing the curves (in pink) selected on Figure 8, we observe two modes of oscillation of the tide: for a fixed "Sea Level" the "Tidal Range" controls the amplitude of oscillation of the tide. The existence of this behavior is coherent with the physics of the problem but it could not be observed in the performed Sobol-based sensitivity analysis.

**SOFTWARE IMPLEMENTATION**

The system described in this article was developed by a collaboration between visualization scientists and statisticians. The aim is the development of mathematical tools to study and analyze multi-run simulations, before integrating the more efficient algorithms in the OpenTURNS software [45]. It was decided to design and implement new interactive visual analytics methods in OpenTURNS by integrating the new developments into ParaView [46]. OpenTURNS and ParaView are both integrated in the SALOME open-source numerical simulation platform [47].

The original idea was to introduce a Functional Boxplot view in ParaView in order to avoid visual clutter and interactively study the outliers of an ensemble of curves. The bivariate PCA plane and the High Density Regions (HDR) were seen as a way of augmenting the information of the Functional Box Plot view, which was an advantage over functional depth methods like [27-28]. Data depth does not allow to display data multimodality but only to calculate quantiles of functions. In our system, if the structure is multi-modal then the analyst can visually identify the clusters, which are disjoint regions of the inner HDR. Furthermore, the bivariate PCA plane could be segmented by any automatic clustering algorithm. This is straightforward in our integration in ParaView because a clustering algorithm can be added to its standard visualization pipeline. We positively tested the Paraview's native implementation of k-means.

On the other side, interacting among views introduced problems in the architecture of ParaView and, as a consequence of this work, the so-called linked views mechanism was developed. In this context, other statistical views were also implemented. We remark the implementation of an interactive Scatter Plot Matrix view that is a version of the work of [41].

Finally, our system was fully integrated in ParaView and is available from version 5.0.1. This software being Open-Source, the examples given in this article can easily be reproduced. Indeed, we include all data presented in this article as supplemental material (https://gitlab.kitware.com/edf/visual-sensitivity-analysis-of-curves).

**CONCLUSION**

We have designed and implemented a system allowing the in-depth study of augmented ensembles issued from multi-run numerical simulations dealing with uncertainty. These augmented ensembles are composed by functions and their associated parameters. The main contribution of our system is that a *visual sensitivity study* becomes possible by jointly analyzing functional outputs and their corresponding input parameters.

Figure 9 synthesizes the overall methodology. Its principal element is based on HDR computed on the PCA bivariate plane. This allows the realization of the following tasks:
- Avoid visual clutter by visualizing interquantile areas and the median curve;
Interactively detect functional outliers;
- Identify clusters of functions by means of the HDR and PCA-plane.

Combining all these elements with the linking of functional outputs to their corresponding input parameters allows the realization of a *visual sensitivity study.*

Two synthetic examples and one industrial use-case have allowed to demonstrate the potential of the approach which has been integrated in a software environment based on the ParaView and OpenTURNS platform. Current works turn to extend this method to larger number of components retained in the PCA step and to the visual sensitivity analysis of parameter-augmented ensembles of spatial fields. Indeed, in a lot of applications, outputs of computer codes are vectors supported by surfaces (see some examples in [2-3, 19-20]). Future works will also consider non-linear dimensionality reduction techniques [31] in order to replace the PCA one.

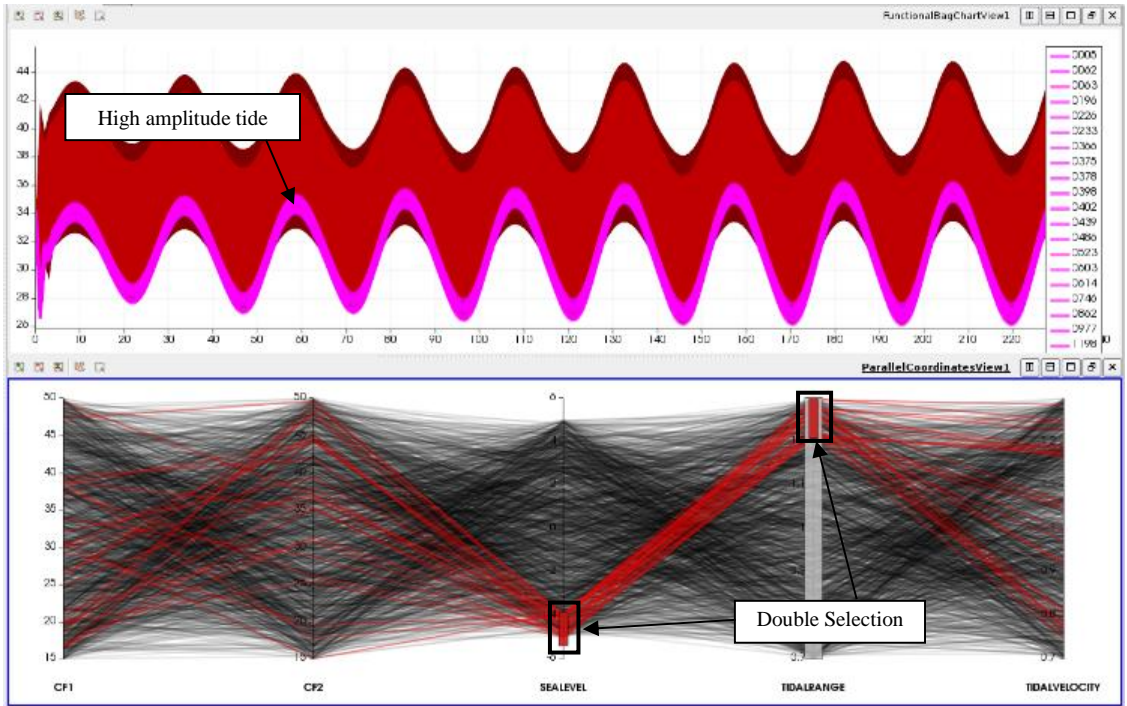

(a)

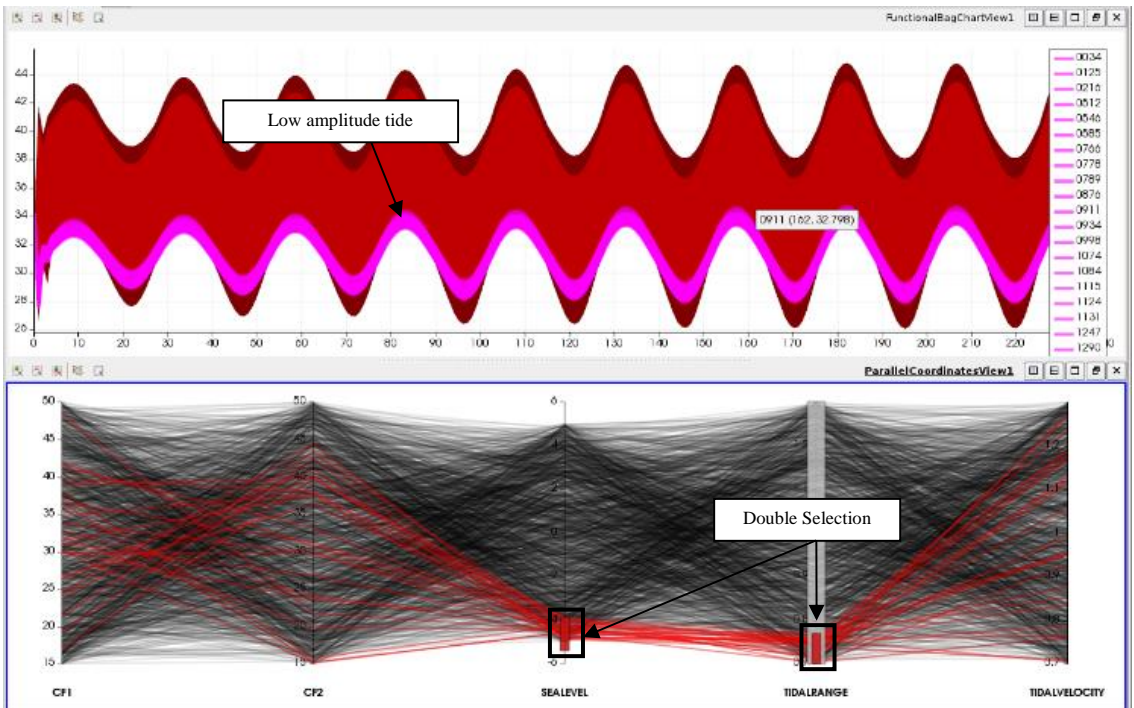

(b)

Figure 8. Interactive exploration of a subtle phenomenon involving two parameters: "Sea Level" and "Tidal Range". The effect of "Tidal Rage" is not the same depending on the value of "Sea Level".

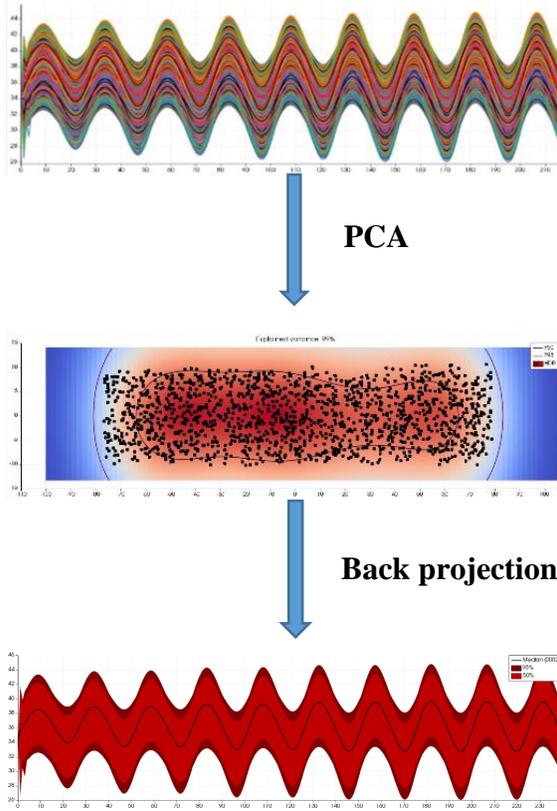

*The ensemble of sampled 1D functions is assembled into a matrix*

**PCA**

-*A kernel-based PDF is calculated on the plane of the two first PCA components*
- *A blue to red colormap is applied to the PDF*
- *HDR are estimated*
- *HDR boundaries are isoprobability contours*

**Back projection**

-*The highest density point estimates the median*
-*HDR boundaries correspond to the functional quantiles*
-*Points outside the outer HDR (95%) gives outliers*

Figure 9. Scheme of the overall PCA-based methodology.